# Two-Step market clearing for local energy trading in feeder-based markets

*Mohsen Khorasany\*, Yateendra Mishra\*, Gerard Ledwich\**

*\*School of Electrical Engineering and Computer Science, Queensland University of Technology, Australia*
*Emails: m.khorasany@qut.edu.au, yateendra.mishra@qut.edu.au, g.ledwich@qut.edu.au*

**Keywords:** Transactive Energy, Local Trading, Market Design, Distributed Optimization.

## Abstract

Recent innovations in Information and Communication Technologies (ICT) provide new opportunities and challenges for integration of distributed energy resources (DERs) into the energy supply system as active market players. By increasing integration of DERs, novel market platform should be designed for these new market players. The designed electricity market should maximize market surplus for consumers and suppliers and provide correct incentives for them to join the market and follow market rules. In this paper, a feeder-based market is proposed for local energy trading among prosumers and consumers in the distribution system. In this market, market players are allowed to share energy with other players in the local market and with neighborhood areas. A Two-Step Market Clearing (2SMC) mechanism is proposed for market clearing, in which in the first step, each local market is cleared independently to determine the market clearing price and in the second step, players can trade energy with neighborhood areas. In comparison to a centralized market, the proposed method is scalable and reduces computation overheads, because instead of clearing market for a large number of players, the market is cleared for a fewer number of players. Also, by applying distributed method and Lagrangian multipliers for market clearing, there is no need for a central computation centre and private information of market players. Case studies demonstrate the efficiency and effectiveness of the proposed market clearing method in increasing social welfare and reducing computation time.

## 1 Introduction

The recent developments of Information and Communication Technologies (ICT) and smart systems has motivated consumers to become more active players instead of being submissive ratepayers. These new active players participate in the market as prosumers by using their local resources, managing their demand, and communicating with other players. Using the two-way flow of information and energy, these new players can exchange both information and energy among themselves and with the grid. In the distribution grid, prosumers can trade energy with other prosumers and consumers by buying energy from them when they need more energy or selling to them when they have surplus energy. This energy trading incentivizes consumers and prosumers to participate more actively in the market. As the number of prosumers with Distributed Energy Resources (DERs) increases in the distribution grid, the reliance on the conventional main electric grid is reduced by providing local generation instead of produced energy in centralized and large energy plants. The local generation of DERs provides opportunities for energy trading in an open and flexible market. This new context of an open and flexible market for DER owners has triggered considerable attention to local energy trading. In the traditional markets, energy trading is performed for large-scale generations and on the transmission side. However, in the local energy market, energy trading occurs between a lot of small-scale DERs. Therefore, applying traditional markets to the local energy market is not straightforward and new mechanisms for small-scale trading and local energy markets are needed. In the local energy markets, energy system is restructured in a way that enables prosumers and consumers to join the energy supply system [1]. A local market is defined based on micro market concept in a residential area or similar for energy trading between people located in the same neighbourhood with the possibility of trading between neighbourhood areas [2].

The design of a market for incentivizing local energy trading between small-scale prosumers and consumers has been investigated in several research. A market design is proposed in [3] to encourage market players to form groups of sellers and buyers for energy trading. The concept of local energy trading between a set of energy storage units has been proposed in [4], where authors propose an auction-based approach for energy trading among storages and consumers. A localized retail energy market is proposed in [5], where an energy broker acts as middleman to match sellers and buyers based on search theory. Authors in [6], propose an auction-based market clearing for local energy trading, where knapsack algorithm is applied for market clearing in P2P energy trading. As more prosumers are integrated into power systems, the centralized methods are not suitable for such heterogeneous systems with large number of players [7] and the scalability and required computation infrastructure would be challenging for any new market design.

In the recent years, there is a tendency for using distributed methods for market design. In these methods, decision making can be distributed across many decision makers which make these methods scalable and reduce computation overheads by eliminating the need for a central computation centre.

Distributed methods need a coordination signal among market players to reach the global optimum, where electricity price or power mismatch can be used as coordination signal. Market design based on distributed methods has gained high attention in the recent literature. Authors in [8] propose a community based electricity market structure which allows prosumers to actively optimize their assets using distributed optimization A distributed consensus-based algorithm for energy trading is proposed in [9], where local estimation of power mismatch is used as coordination signal to reach consensus. Business models for P2P energy trading in distributed energy trading has been proposed in [10], which shows that local energy trading is preferred solution for a system with large number of small-scale DERs. Energy trading among microgrids is investigated in [11] using a distributed solution to minimize social cost of the system. This paper designs a bargaining based energy trading to encourage local energy trading.

In the most of the previous works on distributed energy trading, the focus is on determining the total amount of traded energy to satisfy the overall system demand without paying attention to the computational overheads. However, when the number of market players is high, computational cost and convergence rate become critical factors in the market design and implementation. In this paper, a feeder-based market topology for local energy trading is proposed, where the total grid is divided into several local markets with a limited number of players. In each local market, prosumers and consumers can trade energy to maximize their social welfare. The advantage of the proposed market is that it can be implemented for any number of players by forming any number of local markets that can be cleared locally without any need to individual and private information of players. In other word, instead of clearing the market for a large number of players it can be easily performed for a local market with a limited number of players. It makes the system easier to manage and operate and reduces the burden of communication and computation in the system. Also, a Two-Step Market Clearing (2SMC) approach is presented for market clearing which has the following features:

- There is no need to individual information of market players for market clearing. In addition, it is not required to have a communication channel between consumers and prosumers to exchange information with all players.
- Each player can decide on his action in the market by responding to the coordination signal (price signal)
- Players have the opportunity to trade with neighborhood area to increase their welfare.
- The proposed approach can be implemented for a large number of players.

The rest of this paper is organized as follow: Section 2 gives details of assumption for feeder-based market design and formulation for optimization problem. The 2SMC has been explained in Section 3. Numerical analysis and case studies are provided in Section 4. Finally, Section 5 summarizes the paper and discusses future works.

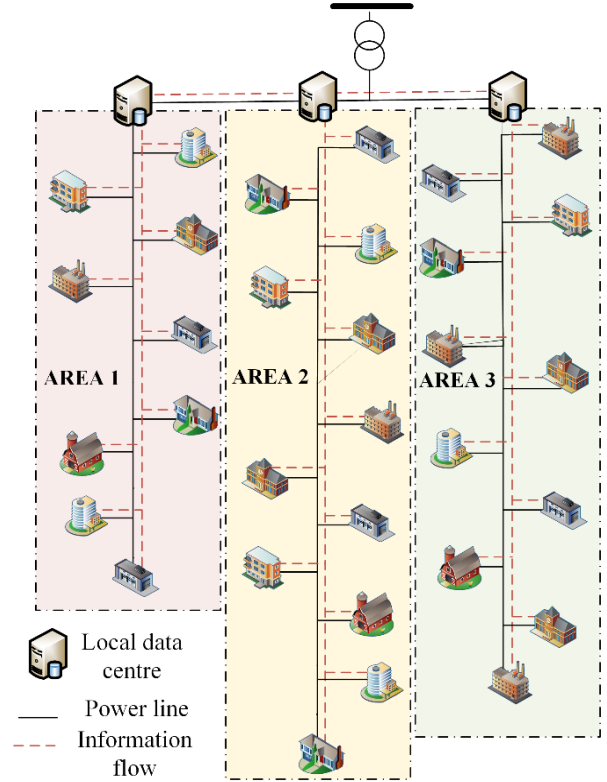

Figure 1: The schematic of feeder-based market for local energy trading

## 2 Feeder-based market

### 2.1 Assumptions

This section provides an overview of the feeder-based energy trading in the distribution network. In this paper, the market structure is designed to incentivize market players for local energy trading. The total distribution system is divided to several local market and prosumures and consumers connected to the same feeder can trade energy locally, where players with excess energy try to sell this energy and players who need energy try to buy it from local market for a lower price. Market players participate in the market by responding to the electricity price and adjusting their load/generation. In each area, there is a data centre which receives information from players of that area to clear the market and acts as coordinator. Figure. 1 illustrates schematic of a feeder-based market with three areas including sellers, buyers and data centres.

### 2.2 Problem formulation

The objective of the local market is to maximize social welfare of all market players. Social welfare maximization not only can maximize the total welfare of all market players, but also can guarantee that each individual welfare is maximized [12]. The number of total areas is specified by $\mathcal{L}$ and $N_\ell^S$ and $N_\ell^B$ indicate the number of sellers (prosumers) and buyers (consumers) in area $\ell$ respectively. The utility function of consumers and cost function of prosumers are defined in this section. The utility function of each consumer is a unique function based on parameters $\omega$ and $\mu$. The utility function for $j^{th}$ consumer in area $\ell$ can be modelled by (1) [13].

$$U^j(d_\ell^j) = \begin{cases} \omega^j d_\ell^j - \mu^j(d_\ell^j)^2 & d_\ell^j < \omega^j/2\mu^j \\ (\omega^j)^2/4\mu^j & d_\ell^j \geq \omega^j/2\mu^j \end{cases} \quad (1)$$

where $d_\ell^j$ is demanded energy by $j^{th}$ consumer from the market. The welfare in buyer side can be measured with the consumer`s surplus which is presented by (2):

$$\mathcal{WB}^j(d_\ell^j) = U^j(d_\ell^j) - \lambda d_\ell^j \quad (2)$$

where $\lambda$ is the paid money by consumer for each kWh of energy. The total cost of output power of prosumers can be estimated by a quadratic cost function with parameters $a$, $b$ and $\gamma$ which customize cost function for each prosumer [14].

$$C^i(s_\ell^i) = a^i(s_\ell^i)^2 + b^i s_\ell^i + \gamma^i \quad (3)$$

The welfare of prosumer is the gained benefit by prosumer in the market and can be presented by (4).

$$\mathcal{WS}^i(s_\ell^i) = \lambda s_\ell^i - C^i(s_\ell^i) \quad (4)$$

where, $s_\ell^i$ is the supplied energy by prosumer $i$ to the market. In each area, prosumers and consumers can trade energy with other players in that area (intra-area) and neighbourhood areas (inter-area). Therefore, the demanded/supplied energy by each consumer/prosumer can be represented by (5) and (6).

$$d_\ell^j = d_{\ell,A}^j + d_{\ell,E}^j \quad (5)$$

$$s_\ell^i = s_{\ell,A}^i + s_{\ell,E}^i \quad (6)$$

where, $d_{\ell,A}^j$ and $d_{\ell,E}^j$ are demanded energy by consumer $j$ in area $\ell$ from intra-area and inter-area markets respectively, whereas $s_{\ell,A}^i$ and $s_{\ell,E}^i$ are the supplied energy by prosumer $i$ in area $\ell$ to these markets respectively. The social welfare in each area is the summation of utility of all consumers minus cost of all prosumers. Therefore, the objective function for all players in the market can be written as (7).

$$\max \sum_{\ell=1}^{\mathcal{L}} \left( \sum_{j=1}^{N_\ell^B} U^j(d_\ell^j) - \sum_{i=1}^{N_\ell^S} C^i(s_\ell^i) \right) \quad (7)$$

Since the demand and supply in the system should be balanced, the problem in (7) is constrained by (8). Also, each player has its own limitation for generation or load as (9) and (10).

$$\sum_{\ell=1}^{\mathcal{L}} \sum_{i=1}^{N_\ell^S} (s_{\ell,A}^i + s_{\ell,E}^i) = \sum_{\ell=1}^{\mathcal{L}} \sum_{j=1}^{N_\ell^B} (d_{\ell,A}^j + d_{\ell,E}^j) \quad (8)$$

$$d_\ell^{j,min} \leq d_{\ell,A}^j + d_{\ell,E}^j \leq d_\ell^{j,max} \quad (9)$$

$$s_\ell^{i,min} \leq s_{\ell,A}^i + s_{\ell,E}^i \leq s_\ell^{i,max} \quad (10)$$

The optimization problem in (1), can be easily solved using different methods such as interior point method or bundle methods. However, this optimization needs individual information of all players. Also, a central computation centre is required to collect information of all players and find energy allocation and price for all of them. With the increase in the number of market players, the required computation facilities and communication signals for the market clearing would a barrier for new market design. Therefore, in this paper, a distributed method is applied for market clearing, which can be implemented without any need to central computation centre for all players. In this approach, the total problem is decomposed to several sub-problems for each area, and each sub-problem is also decomposed to several local problems which are locally solvable by each player.

## 3 Two-step market clearing approach

In this paper, distributed market clearing approach utilizing the primal-dual gradient descent method is used for market clearing [15]. In this approach, the objective function is augmented with Karush-Kuhn-Tucker (KKT) multipliers and after that, a distributed iterative approach is developed that moves in the direction of maximizing the augmented objective function [16, 17]. Here, the KKT multipliers are coordination signals which should be updated iteratively till are players reach consensus on their actions in the market. The Lagrangian function to solve (7) subject to (8) can be defined as follows:

$$\begin{aligned} &L(s_{\ell,A}^i, s_{\ell,E}^i, d_{\ell,A}^j, d_{\ell,E}^j) = \\ &\sum_{\ell=1}^{\mathcal{L}} \left( \sum_{j=1}^{N_\ell^B} U^j(d_\ell^j) - \sum_{i=1}^{N_\ell^S} C^i(s_\ell^i) \right) \\ &+ \lambda_\ell \left( \sum_{\ell=1}^{\mathcal{L}} \sum_{i=1}^{N_\ell^S} s_{\ell,A}^i - \sum_{\ell=1}^{\mathcal{L}} \sum_{j=1}^{N_\ell^B} d_{\ell,A}^j \right) \\ &+ \lambda_C \left( \sum_{\ell=1}^{\mathcal{L}} \sum_{i=1}^{N_\ell^S} s_{\ell,E}^i - \sum_{\ell=1}^{\mathcal{L}} \sum_{j=1}^{N_\ell^B} d_{\ell,E}^j \right) \end{aligned} \quad (11)$$

where, $\lambda_\ell$ and $\lambda_C$ are KKT multipliers. The next step is to develop the distributed iterative approach based on dual decomposition. In this market, each prosumer/consumer can sell/buy energy to/from local market or neighborhood markets. Market clearing is designed to be performed in two step. At the first step, players can trade energy intra-area with other players with a unique clearing price for each area ($\lambda_\ell$). In the second step, players can trade energy with neighborhood areas, where prosumers of an area with lower price can sell their energy to the consumers of an area with higher price. This two-step market clearing allows market players to decide on their action in the market by adjusting their load/generation in response to the clearing price. Also, it gives them opportunity to increase their welfare through trading with neighbourhood areas. Therefore, the iterative approach for market clearing should be developed to clear the market in two step. At the first step, prosumers and consumers in each area send their demand and supply to the local data centre. The electricity price in each area is calculated based on all demands and supplies using (12).

$$\lambda_\ell(k+1) = \lambda_\ell(k) + \eta_\ell \left( \sum_{i=1}^{N_\ell^S} s_{\ell,A}^i(k) - \sum_{j=1}^{N_\ell^B} d_{\ell,A}^j(k) \right) \quad (12)$$

where $\eta_\ell$ is tuning parameter which indicates convergence rate of demand and supply and $\lambda_\ell$ is clearing price for intra-area trading. Then, this price is sent to the players and they will update their demand and supply based on this price by (13) and (14) respectively.

$$s_{\ell,A}^i = \max_{s_\ell^{i,min} \leq s_{\ell,A}^i \leq s_\ell^{i,max}} \left(s_{\ell,A}^i(\lambda_\ell(k+1)) - C^i(s_{\ell,A}^i)\right) \quad (13)$$

$$d_{\ell,A}^j = \max_{d_\ell^{j,min} \leq d_{\ell,A}^j \leq d_\ell^{j,max}} \left(U^j(d_{\ell,A}^j) - d_{\ell,A}^j(\lambda_\ell(k+1))\right) \quad (14)$$

This algorithm repeats till demand and supply are converged. In the second step of market clearing, a new area (indicated by C) is formed by prosumers of the area with lower price and consumers of the area with higher price and the described procedure for the first step will repeat for the market clearing in the new area using (15-17).

$$\lambda_C(k+1) = \lambda_C(k) + \eta_C \left(\sum_{i=1}^{N_C^S} s_{\ell,E}^i(k) - \sum_{j=1}^{N_C^B} d_{\ell,E}^j(k)\right) \quad (15)$$

$$s_{\ell,E}^i = \max_{s_\ell^{i,min} \leq s_{\ell,A}^i + s_{\ell,E}^i \leq s_\ell^{i,max}} \left(s_{\ell,E}^i(\lambda_C(k+1)) - C^i(s_{\ell,E}^i)\right) \quad (16)$$

$$d_{\ell,E}^j = \max_{d_\ell^{j,min} \leq d_{\ell,A}^j + d_{\ell,E}^j \leq d_\ell^{j,max}} \left(U^j(d_{\ell,E}^j) - d_{\ell,E}^j(\lambda_C(k+1))\right) \quad (17)$$

where, $\eta_C$ is tuning parameter for the second step and $\lambda_C$ is clearing price for inter-area trading. In the second step, a new market clearing price in calculated using (15) and sellers and buyers respond to this price by managing their generation and load using (16) and (17) respectively. Market clearing algorithm is shown in **Algorithm 1**. In this algorithm, $\lambda_\ell$ and $\lambda_C$ are actually Lagrangian multiplier related to demand-supply constraint which act as coordination signals for distributed market clearing. These prices are calculated in data centres and any calculation for updating demand and supply is performed in market players' smart meters.

---

**Algorithm 1: Market clearing algorithm**

**Initialization:** Set $\lambda_\ell(0)$, $\epsilon_\ell$ and $\eta_\ell$

*Step 1 (Intra-Area clearing)*
1: **Repeat** till $|\lambda_\ell(k+1) - \lambda_\ell(k)| \leq \epsilon_\ell$
2: Update $\lambda_\ell$ according to
$$\lambda_\ell(k+1) = \lambda_\ell(k) + \eta_\ell \left(\sum_{i=1}^{N_\ell^S} s_{\ell,A}^i(k) - \sum_{j=1}^{N_\ell^B} d_{\ell,A}^j(k)\right)$$
3: Send $\lambda_\ell(k+1)$ to the players in the area
4: Sellers update their supply according to
$$s_{\ell,A}^i = \max_{s_\ell^{i,min} \leq s_{\ell,A}^i \leq s_\ell^{i,max}} \left(s_{\ell,A}^i(\lambda_\ell(k+1)) - C^i(s_{\ell,A}^i)\right)$$
5: Buyers update their demand according to
$$d_{\ell,A}^j = \max_{d_\ell^{j,min} \leq d_{\ell,A}^j \leq d_\ell^{j,max}} \left(U^j(d_{\ell,A}^j) - d_{\ell,A}^j(\lambda_\ell(k+1))\right)$$

*Step 2 (Inter-Area clearing)*
1: Communicate with neighborhood areas
2: Set area C with sellers of area with lower price and buyers of area with higher price
3: Set $\lambda_C(0)$, $\epsilon_C$ and $\eta_C$
4: **Repeat** till $|\lambda_C(k+1) - \lambda_C(k)| \leq \epsilon_C$
5: Update $\lambda_C$ according to
$$\lambda_C(k+1) = \lambda_C(k) + \eta_C \left(\sum_{i=1}^{N_C^S} s_{\ell,E}^i(k) - \sum_{j=1}^{N_C^B} d_{\ell,E}^j(k)\right)$$
6: Send $\lambda_C(k+1)$ to the players in the area
7: Sellers update their supply according to
$$s_{\ell,E}^i = \max_{s_\ell^{i,min} \leq s_{\ell,A}^i + s_{\ell,E}^i \leq s_\ell^{i,max}} \left(s_{\ell,E}^i(\lambda_C(k+1)) - C^i(s_{\ell,E}^i)\right)$$
8: Buyers update their demand according to
$$d_{\ell,E}^j = \max_{d_\ell^{j,min} \leq d_{\ell,A}^j + d_{\ell,E}^j \leq d_\ell^{j,max}} \left(U^j(d_{\ell,E}^j) - d_{\ell,E}^j(\lambda_C(k+1))\right)$$

---

## 4 Case studies

The proposed 2SMC is implemented in a simple market with 20 players (9 prosumers and 11 consumers) which are dispatched in three areas. In each area, there are a different number of prosumers and consumers that try to maximize their welfare by adjusting their generation/load in response to the market price. Table 1 shows the parameters of the cost function and utility function of consumers and prosumers [18]. The initial value of $\lambda_\ell$ and $\lambda_C$ are set to zero. Also, the minimum demand and supply for all players are selected to be zero. The evolution of the electricity price in step 1 for different areas is shown in Figure.2. This figure illustrates that the clearing price in area 2 ($\lambda_2$) is higher than the price in area 1 ($\lambda_1$) and 3 ($\lambda_3$). Therefore, consumers of area 2 and prosumers of area 1 and 3 participate in the step 2 clearing. Market clearing price for step 2 ($\lambda_C$) is also shown in Figure 2 which is higher than $\lambda_1$ and $\lambda_3$ and lower than $\lambda_2$.

The proposed method is compared with one-step market clearing (1SMC) method, where there is only one market clearing price for all market players. In this case, players of different areas participate in the one market and market clearing is performed in one step. Market clearing price for 1SMC is shown in Figure 1 as $\lambda_T$. Figure 3 shows the social welfare of market players and traded energy for 1SMC and 2SMC. These results reveal that by applying the proposed 2SMC, total traded energy and social welfare of market players are increased, which means that in the proposed market clearing, players gain more welfare in the market.

In order to analyze performance of the proposed algorithm for a large scale system, it has been applied to a system with 2000 players, including 900 prosumers and 1100 consumers. For this system, the computational time and social welfare of the proposed 2SMC are compared with the 1SMC and results are shown in Table 2. In 2SMC, the computational time is obtained by adding the highest time in the first step and the required time in the second step.

| Prosumers` parameters | | | | | Consumers` parameters | | | | |
|---|---|---|---|---|---|---|---|---|---|
| $i$ | $a^i$ | $b^i$ | $\ell$ | $s_\ell^{i,max}$ | $j$ | $\omega^j$ | $\mu^j$ | $\ell$ | $d_\ell^{j,max}$ |
| 1 | 0.0031 | 8.71 | 1 | 113.23 | 1 | 17.17 | 0.0935 | 1 | 91.79 |
| 2 | 0.0074 | 3.53 | 1 | 179.1 | 2 | 12.28 | 0.0417 | 1 | 147.29 |
| 3 | 0.0066 | 7.58 | 1 | 90.03 | 3 | 18.42 | 0.1007 | 1 | 91.41 |
| 4 | 0.0063 | 2.24 | 2 | 106.41 | 4 | 7.06 | 0.0561 | 1 | 62.96 |
| 5 | 0.0069 | 8.53 | 2 | 193.80 | 5 | 10.85 | 0.0540 | 2 | 100.53 |
| 6 | 0.0014 | 2.25 | 3 | 37.19 | 6 | 18.91 | 0.1414 | 2 | 66.88 |
| 7 | 0.0041 | 6.29 | 3 | 195.4 | 7 | 18.76 | 0.0793 | 2 | 118.35 |
| 8 | 0.0051 | 4.30 | 3 | 62.17 | 8 | 15.70 | 0.1064 | 2 | 73.81 |
| 9 | 0.0032 | 8.26 | 3 | 143.41 | 9 | 14.28 | 0.0580 | 2 | 84.00 |
| | | | | | 10 | 10.15 | 0.0460 | 3 | 110.32 |
| | | | | | 11 | 19.04 | 0.0650 | 3 | 146.46 |

Table 1: Utility and cost function parameters of prosumers and consumers

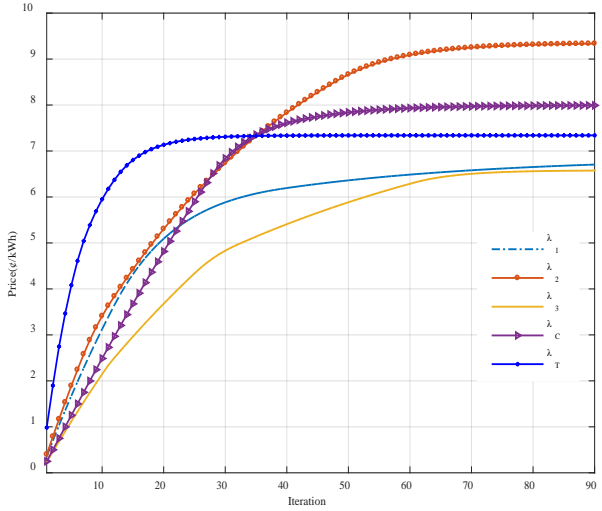

Figure 2: Market clearing prices

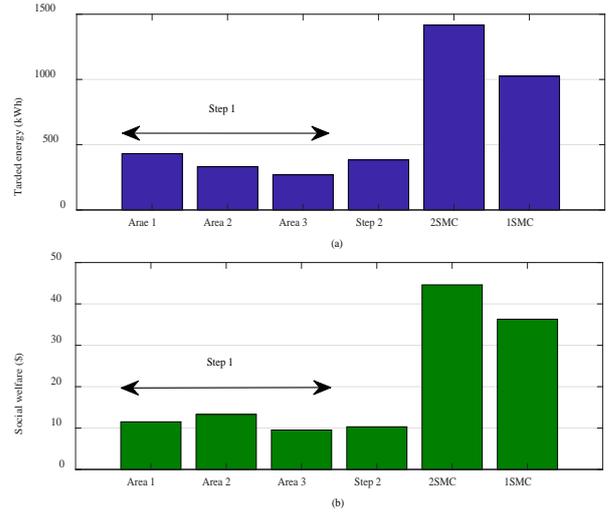

Figure 3: Comparison of 1SMC and 2SMC; (a) Traded energy (b) Social welfare

|  | 2SMC | 1SMC | Percentage of variation (%) |
|---|---|---|---|
| Social Welfare ($) | 4364.2 | 3629.8 | + 20% |
| Computational Time (s) | 19.06 | 37.13 | - 48% |

Table 2: Comparison of 2SMC and 1SMC performance in large-scale system

These results show that for a large scale system, the proposed method can increase the social welfare by 20% and reduce the required computational time by 48%, which validates the efficacy of the proposed method in the large-scale systems.

## 5 Conclusions

This paper proposes a feeder-based market topology for local energy trading among prosumers and consumers. In the feeder-based markets, the total grid is divided into several local markets with limited number of players and in each local market, prosumers and consumers can trade energy to maximize their social welfare. The advantage of the proposed market is that it can be implemented for any number of players by forming any number of local markets that can be cleared locally without any need to individual and private information of players. Then a new approach for market clearing in feeder based market is proposed, where a 2SMC mechanism is designed to clear the market for local energy trading. This approach can be performed without any need to central computation centre and individual information of players. Simulation results verified that the proposed 2SMC increases total traded energy and welfare of market players. Also, this method reduces the computational time and can be implemented for any large-scale system. For the future works, the networks constraints (especially line flow constraint) will be added for trading among neighbouring areas to consider network topology in the market clearing. Also, other market players such as retailers and utility companies will be modelled in the market.